\documentclass[12pt]{article}
\usepackage{amsmath, amsthm, amssymb}
\usepackage{amsfonts}
\usepackage{epsfig}
\usepackage{graphics}
\usepackage{subfigure}
\usepackage{graphicx}
\usepackage{color}
\usepackage{epstopdf} 
\usepackage{float} 
\usepackage{mathrsfs} 
\usepackage{relsize} 
\newcommand{\ds}{\displaystyle}

\newcommand{\noi}{\noindent}
\newcommand{\bq}{\begin{equation}}
\newcommand{\eq}{\end{equation}}
\newcommand{\bqr}{\begin{eqnarray}}
\newcommand{\eqr}{\end{eqnarray}}
\newcommand{\bqrn}{\begin{eqnarray*}}
\newcommand{\eqrn}{\end{eqnarray*}}

\begin{document}

\begin{center}	
{\bf Dual Approach as Empirical Reliability for Fractional Differential Equations}\\[2mm]
 Pavel B. Dubovski, Jeffrey Slepoi\\
The paper has been published in \emph{Journal of Physics: Conference Series {\bf 2099}} (2021); DOI:10.1088/174206596/2099/1/012004
\end{center}

\begin{abstract}
	Computational methods for fractional differential equations exhibit essential instability.  Even a minor modification of the coefficients or other entry data may switch good results to the divergent.  The goal of this paper is to suggest the reliable dual approach which fixes this inconsistency.  We suggest to use two parallel methods based on the transformation of fractional derivatives through integration by parts or by means of  substitution. We introduce the method of substitution and choose the proper discretization scheme that fits the grid points for the by-parts method. The solution is reliable only if both methods produce the same results. As an additional control tool, the Taylor series expansion allows to estimate the approximation errors for fractional derivatives. In order to demonstrate the proposed dual approach, we apply it to linear, quasilinear and semilinear equations and obtain very good precision of the results. The provided examples and counterexamples support the necessity to use the dual approach because either method, used separately, may produce incorrect results. The order of the exactness is close to the exactness of fractional derivatives approximations.
\end{abstract}

\noindent Keywords: fractional differential equations; empirical reliability; substitution method; dual approach 

\section{Introduction} 

There are several definitions of fractional derivatives, of which the best known are the Riemann-Liouville and Caputo definitions.  Caputo's definition (1967) of fractional derivative of order $\alpha>0$ is especially well-suited for initial value problems:\\
\begin{equation}\label{MainFD}
 D^\alpha f(t) := \frac{1}{\Gamma(n-\alpha)}\int_{0}^{t} \frac{f^{(n)}(x)}{(t-x)^{\alpha-n+1}}dx, 
\end{equation}
where $n-1 \le \alpha < n, n \in \mathbb{N}$ and $\Gamma(n-\alpha)$ is the gamma function.

It is worth mentioning that Caputo derivative was first derived by Soviet mechanic A.N.Gerasimov in 1948 \cite{Gerasimov}. 

The main obstacle to numerical solving FrDE is the presence of the mild singularity in its definition. Integrating the fractional derivative by parts helps to fix this problem. This method (by-parts method) has been used, e.g., in \cite{Ramzi} and \cite{Ramzi2}, where the Caputo fractional derivative was approximated by finite difference scheme. Many observations demonstrate that fractional derivatives can be calculated with high precision.
For benchmarking tests, like in \cite{Ramzi} and \cite{Ramzi2}, the numerical results exhibit good convergence to exact solutions. However, if the equations are slightly modified, the convergence fails quite often as we show in Section \ref{SectFailedExamples}. In these cases, we cannot guarantee the convergence to the correct solution. In Figures 1b and 2 
we show simple examples of such a divergence, which may happen even for linear equations.

The goal of this paper is to develop a reliable dual approach to find a numerical solution to fractional differential equations. 
As sufficiently general samples, we consider quasilinear fractional differential equations of the form
\begin{equation}\label{QLEquation}
\sum_{k=i}^{n}K_i(u,x)D^{\alpha_i} u(x) + f(x)=g(u(x)).
\end{equation}
During discretization we utilize forward,  backward, and central differences with higher precision \cite{Mathews}.

To address the convergence issues, we apply two methods of resolving the aforementioned singularity: the method of substitution and the by-parts method. Both methods produce almost-identical correct results when the solution is stable and, especially, is known analytically.  However, we show in Section \ref{SectFailedExamples} that after seemingly minor modifications to the equation, either one or both methods fail to match the analytical solutions. Our analysis leads us to conclude that the computations are reliable if both methods produce the same result. For known functions, our analysis is supported by using the Taylor series.

We demonstrate how the dual approach can be used to find solutions and evaluate their reliability for linear and nonlinear fractional differential equations \eqref{QLEquation}.

In Section \ref{SectSubst} we introduce the substitution method, in Section \ref{SectByParts} we describe the well-known method of integration by parts (by-parts method), in Section \ref{SectDescrit} we describe the methods of discretization.  Next Section \ref{SectFailedExamples} is critical: we provide the examples of linear equations with the by-parts and substitution methods. We demonstrate that only if the results of the methods match they yield the reliable solution. As an example of the insufficiency of just one method in the computations of fractional differential equations, we reconsider an example from the well-known monograph [\cite{Podlubny}, Section 8.4, Example 3] and make the necessary corrections. Section \ref{SectExamples} contains the examples of quasilinear fractional equations with essential nonlinearity in the right-hand side. 

\section {Method of substitution}\label{SectSubst}

Let us convert the Caputo fractional derivative into a simpler form by substituting in \eqref{MainFD}\\ $\alpha-n+1=\frac{b}{b+1}$ and $u=(t-x)^{\frac{1}{b+1}}$, where $n-1<\alpha<n$ and therefore $0<\alpha-n+1<1$:
\begin{equation*}
b= \frac{\alpha-n+1}{n-\alpha} > 0
\end{equation*}
This yields:
\begin{equation}
	du = -\frac{1}{b+1}(t-x)^{-\frac{b}{b+1}}dx=\frac{1}{b+1}\cdot \frac{1}{(t-x)^{\alpha-n+1}}dx; \text { }
f^{(n)}(x)=f^{(n)}(t-u^{b+1});
\end{equation}
and therefore
\begin{eqnarray}\label{NoSingForm}
D^\alpha f(t)&  = &\frac{1}{\Gamma(n-\alpha)}\int_{0}^{t} \frac{f^{(n)}(x)}{(t-x)^{\alpha-n+1}}dx\nonumber \\
 & = &-\frac{b+1}{\Gamma(n-\alpha)}
\int_{t^{\frac{1}{b+1}}}^{0}f^{(n)}(t-u^{b+1})du\nonumber\\
& = &\frac{1}{\Gamma(n-\alpha)}\cdot\frac{1}{n-\alpha} \int_{0}^{t^{n-\alpha}}f^{(n)}(t-u^{\frac{1}{n-\alpha}})du.
\end{eqnarray}
In the last step we used $b+1=(n-\alpha)^{-1}$. 

The first natural wish for numerical integration is splitting the interval $0<u<t^{n-\alpha}$ into equal parts. 
However, for each fractional derivative the step becomes different (since it depends on $\alpha$). Hence, if the equation contains different derivatives, we have to deal with different grid points for each derivative, and this way is not applicable for equations with several fractional derivatives.

To fix this difficulty, we calculate $f^{(n)}(x_k)$ at each point $x_k=kh$ between $0$ and $t=mh$ with step $h$
and
introduce $u_k=(t-x_k)^{n-\alpha}$. 
Then the grid points $x_k=kh$ fit to all derivatives and functions in the equation, and the use of the trapezoid rule yields  
\begin{equation}\label{FrDSubst}
D^\alpha f(t) \approx \frac{1}{\Gamma(n+1-\alpha)} \sum_{k=1}^{m}\frac{f^{(n)}(x_k)+f^{(n)}(x_{k-1})}{2}\cdot(u_{k-1}-u_{k}).
\end{equation}
%
Since the increase in $x$ leads to the decrease in $u$, we obtain $u_{k-1}>u_k$.\\
This approach implies that Caputo derivative can be implemented as follows:
\begin{equation}\label{GL-Caputo-def} 
D^\alpha f(t) = \frac{1}{\Gamma(n+1-\alpha)} \lim\limits_{m\to \infty}  \sum_{k=1}^{m}\frac{f^{(n)}(x_k)+f^{(n)}(x_{k-1})}{2}\left ((t-x_{k-1})^{n-\alpha}-(t-x_k)^{n-\alpha}\right),
\end{equation}
where $n=\lfloor \alpha \rfloor + 1$ and, hence, $n-1 \le \alpha < n,\ n \in \mathbb{N}$.

\section{Integration by parts}\label{SectByParts} 
In this well-known method the singularity is eliminated by integration by parts
\bqr \label{FDApprox} 
D^\alpha f(t) &=& \frac{1}{\Gamma(n-\alpha)}\int_0^{t} \frac{f^{(n)}(x)}{(t-x)^{\alpha-n+1}}dx \nonumber\\
   &=& \frac{1}{\Gamma(n+1-\alpha)}\left(f^{(n)}(0)t^{n-\alpha}+\int_0^{t} (t-x)^{n-\alpha}f^{(n+1)}(x)dx\right),
\eqr
where $n-1< \alpha < n$ and after the elimination of singularity, the integral can be approximated by using the trapezoidal rule as follows
\begin{equation}\label{FDIntApprox}
\int_0^{t} (t-x)^{n-\alpha}f^{(n+1)}(x))dx \approx \frac{h}{2}\left[ t^{n-\alpha}f^{(n+1)}(0) + 2\sum_{j=1}^{k}(t-x_j)^{n-\alpha}f^{(n+1)}(x_j) \right],
\end{equation}
where $x_k=t-h$ since at $x_k=t$ the last term of the trapezoidal rule disappears. 
Finally, in place of expression (\ref{GL-Caputo-def}), for the by-parts method we obtain
\begin{equation} \label{FrDByParts}
D^\alpha f(t) = \frac{1}{\Gamma(n+1-\alpha)}\lim\limits_{h\to 0}\frac{h}{2}\left[ t^{n-\alpha}f^{(n+1)}(0) + 2\sum_{j=1}^{k}(t-x_j)^{n-\alpha}f^{(n+1)}(x_j) \right].
\end{equation}
As we can see from \eqref{GL-Caputo-def} and \eqref{FrDByParts}, the principal difference in the approximations of fractional derivatives is the use of $(n+1)$-st derivative in the by-parts method whereas the highest derivatives used in the substitution method, is of the $n$-th order.
\section{Discretization methods} \label{SectDescrit} 

\noi {\bf 4.1 Numerical discretization for the substitution method}\\
If equation \eqref{QLEquation} is of the first order with $0<\alpha<1$, then we need to represent first derivatives, apply \eqref{FrDSubst}, and utilize the implicit method by setting up a system of algebraic equations. We can calibrate the number of equations for the best performance but the basic idea is the same. Our algorithm provides the second order of accuracy. The steps are as follows: 
\begin{enumerate}
	\item use finite differences to represent a few consecutive derivatives: 
		\begin{subequations}
			\begin{align}\label{DiscrSubst}	
			u'_{k-1} & =\frac{-3u_{k-1}+4u_k-u_{k+1}}{2h}+O(h^2) \text{ -- forward difference}, \\
			u'_k & 	   =\frac{-u_{k+2}+8u_{k+1}-8u_{k-1}+u_{k-2}}{12h}+O(h^2) \label{DS2} \text { -- central difference},\\  	
			u'_{k+1} &  =\frac{3u_{k+1}-4u_{k}+u_{k-1}}{2h}+O(h^2) \text { -- backward difference}. \label{DS3}
			\end{align}
		\end{subequations}
	\item using formula \eqref{FrDSubst} represent fractional derivative $D^{\alpha_i} u(x)$;
	\item multiply each derivative by the corresponding $K_i(x_k,u_k)$;
	\item sum up all expressions $K_i(x_k,u_k)D^{\alpha_i} u_k(x_k)$ and
	add to the above generated expression the value of $f(x_k)$  for each point;
	\item equate the produced sum to $g(u_k)$ for each point.
\end{enumerate}

This process generates the system of algebraic equations. 
It is important to point out that for the equations of the second order with at least one $\alpha_i \in (1,2)$, in addition to the initial condition $u_0=u(0)$, the second initial condition is necessary $u'_0=u'(0)$, which we code using forward approximation $u'_0=\ds\frac{-3u_0+4u_1-u_2}{2h}$. In this case the second derivatives need to be represented like in \eqref{DP1}-\eqref{DP3} (in addition to \eqref{DiscrSubst}-\eqref{DS3}). For higher derivatives their representations can be found, say, in \cite{Mathews}.\\

\noi {\bf 4.2 Numerical discretization for the by-parts method}\\
The steps to solve equation like \eqref{QLEquation} using the by-parts method, are the same as in the substitution method, but the fractional derivative representation changes. For $0<\alpha<1$ first and second derivatives of the unknown function $u(x)$ need to be used. For equations \eqref{FDApprox} and \eqref{FDIntApprox} FDM utilized are (points prior to point $u_k$ are assumed to be already found):
	\begin{subequations}\label{DiscrByParts}
	\begin{align}	
	u'_0 & =\frac{-3u_{0}+4u_{1}-u_{2}}{2h} +O(h^2) \text { -- forward difference}, \\  	
	u''_0 &  =\frac{2u_{0}-5u_{1}+4u_{2}-u_{3}}{h^2}+O(h^2) \text { -- forward difference}, \label{DP1}\\
	&\text {Remaining derivatives in the sum of equation \eqref{FDIntApprox}:}\nonumber\\
	u''_{m} & = \frac{-u_{m-2}+16u_{m-1}-30u_{m}+16u_{m+1}-u_{m+2}}{12h^2}+O(h^4) \text { -- central differences}, \label{DP2} \\
	&\text{ for $m =k-1,k,k+1,k+2$},\nonumber\\
	u''_{k+3} & = \frac{2u_{k+3}-5u_{k+2}+4u_{k+1}-u_k}{h^2}+O(h^2) \text { -- backward difference}.  \label{DP3}
	\end{align}
\end{subequations}
The use of these derivatives in formulas \eqref{FDApprox} and \eqref{FDIntApprox} allows us to represent fractional derivatives in quasilinear equation \eqref{QLEquation}.  Then we apply the steps like in substitution method aboveand generate the algebraic system of equations for $u_k$.  The repetition of this process through the whole interval provides the values of function $u(x)$ on the entire interval.\\

\noi {\bf  4.3 Approximation errors}\\
The above mentioned methods of discretization for both substitution and by parts methods lead to the same level of accuracy. 
In fact, the best way to estimate the precision of calculation of fractional derivative is to compare it with the exact expression.  If it is not available, then the Taylor series approximation can be used when the solution function is known.  In this case, the Caputo derivative for the function $f(x)=\ds\sum_{n=0}^{\infty}f^{(n)}(0)\frac{x^n}{n!}$ can be expressed as follows
\begin{equation}
D^{\alpha} f(t) = \frac{1}{\Gamma(n-\alpha)} \sum_{k>\alpha} \int_{0}^{t} \frac{f^{(k)}(0)x^{k-1}}{(k-1)!(t-x)^{\alpha-n+1}}dx.
\end{equation} 
Then, assuming the convergence of series, we obtain
\begin{equation}                                                                        \label{FDasRowFinal}
D^{\alpha}f(x) = \sum_{k>\alpha}\frac{f^{(k)}(0)\cdot x^{k-\alpha}}{\Gamma(k+1-\alpha)}. 
\end{equation}
Although we could not find any direct reference to formula \eqref{FDasRowFinal}, we believe that it is known and is a kind of a "folklore" in the mathematical community. \\
Now, we are ready to check the accuracy of both methods to calculate the fractional derivative.  For example, Table \ref{TanCompTbl} presents the accuracy of calculations $D^{0.4}\tan x$:  
\begin{table}[H]
	\begin{center}
		\begin{tabular}{ | c | r | r | c | r | c |}\hline
		$x$ & Taylor exp & Substitution & Abs Err Subst & Int. by Parts & Abs Err by Parts   \\ \hline
			0.1 	 & 	 0.2824821555 	 & 	 0.2824821407 	 & 	$1.5\times 10^{-8}$	 & 	 0.2824821402 	 & 	$1.5\times 10^{-8}$	\\\hline
			0.2 	 & 	 0.4344599870 	 & 	 0.4344599557 	 & 	$3.1\times 10^{-8}$	 & 	 0.4344599549 	 & 	$3.2\times 10^{-8}$	\\\hline
			0.3 	 & 	 0.5680457063 	 & 	 0.5680456557 	 & 	$5.1\times 10^{-8}$	 & 	 0.5680456546 	 & 	$5.2\times 10^{-8}$	\\\hline
			0.4 	 & 	 0.6996788619 	 & 	 0.6996787873 	 & 	$7.5\times 10^{-8}$	 & 	 0.6996787858 	 & 	$7.6\times 10^{-8}$	\\\hline
			0.5 	 & 	 0.8392329447 	 & 	 0.8392328384 	 & 	$1.1\times 10^{-7}$	 & 	 0.8392328364 	 & 	$1.1\times 10^{-7}$	\\\hline
			0.6 	 & 	 0.9959906149 	 & 	 0.9959904642 	 & 	$1.5\times 10^{-7}$	 & 	 0.9959904614 	 & 	$1.5\times 10^{-7}$	\\\hline
		\end{tabular}
		\caption{Calculation of fractional derivatives $D^{0.4}\tan x$ with step $h=0.0001$.} \label{TanCompTbl}
	\end{center}	
\end{table}
\noi It is easy to see that we have the accuracy of order $O(h^2)$ as expected.  Similar accuracy was observed for other smooth functions. 
Of course, usually we do not know the Taylor series representation of the solution to differential equations. In this case we can calculate the difference between solutions as a proxy of the error and compare it with the expected accuracy of the calculation of the fractional derivatives (in the presented schemes the accuracy is $O(h^2)$). We resulting order of accuracy of FrDE cannot be expected better than the accuracy of the calculation of the fractional derivatives.

The dual approach works even if the equation contains integer derivatives. In this case, integer derivatives can be approximated by tiny deviation $\alpha=n-\delta$, where $\delta\sim 10^{-14}$, and the accuracy of such the deviation of the order of derivative is almost perfect.

\section{Valid and failed convergences for linear equations}\label{SectFailedExamples} 
The use of only one numerical method to solve fractional equations is not sufficient.  Having two significantly different methods allows us to compare their outputs and, if they match, conclude that the found solution is reliable. 
Four examples, reflected in Figures 1 and 2, 
correspond to linear equation with exact solution $u=x^\beta$
\bq
K(x)D^{0.5}u(x)-u(x)+x^\beta - K(x)D^{0.5}x^\beta=0. \label{5-0}
\eq
In Figure 1, 
$K(x)\equiv 1$, $u(x)=x^{1.2}$ (left) and $u(x)=x^{0.5}$ (right).
In case 1a equation (\ref{5-0}) looks as follows:
\bq
 D^{0.5}u(x)-u(x)+x^{1.2}-1.2\mathcal{B}(0.5,1.2)\frac {x^{0.7}}{\sqrt\pi}=0 \label{5-1a}
\eq
and has exact solution $u(x)=x^{1.2}$. Here $\mathcal{B}$ is beta-function. In case 1b, equation (\ref{5-0}) becomes
\bq
D^{0.5}u(x)-u(x)+x^{0.5}-\frac{\sqrt\pi}{2}=0 \label{5-1b}
\eq
and has exact solution $u(x)=x^{0.5}$.

\noi In Figure 1a the results for both substitution and by parts methods coincide and provide the correct solution $u(x)=x^{1.2}$. The consilience of these results justifies the proposed dual approach that guarantees the reliability of computations.\\
In Figure 1b the results for both substitution and by parts methods differ from the exact solution 
$u(x)=x^{0.5}$ and, most important, are far from each other. Thus, the computations are not reliable.
 
\begin{figure}[H] 
	\begin{center}
		\includegraphics[width=7.5cm]{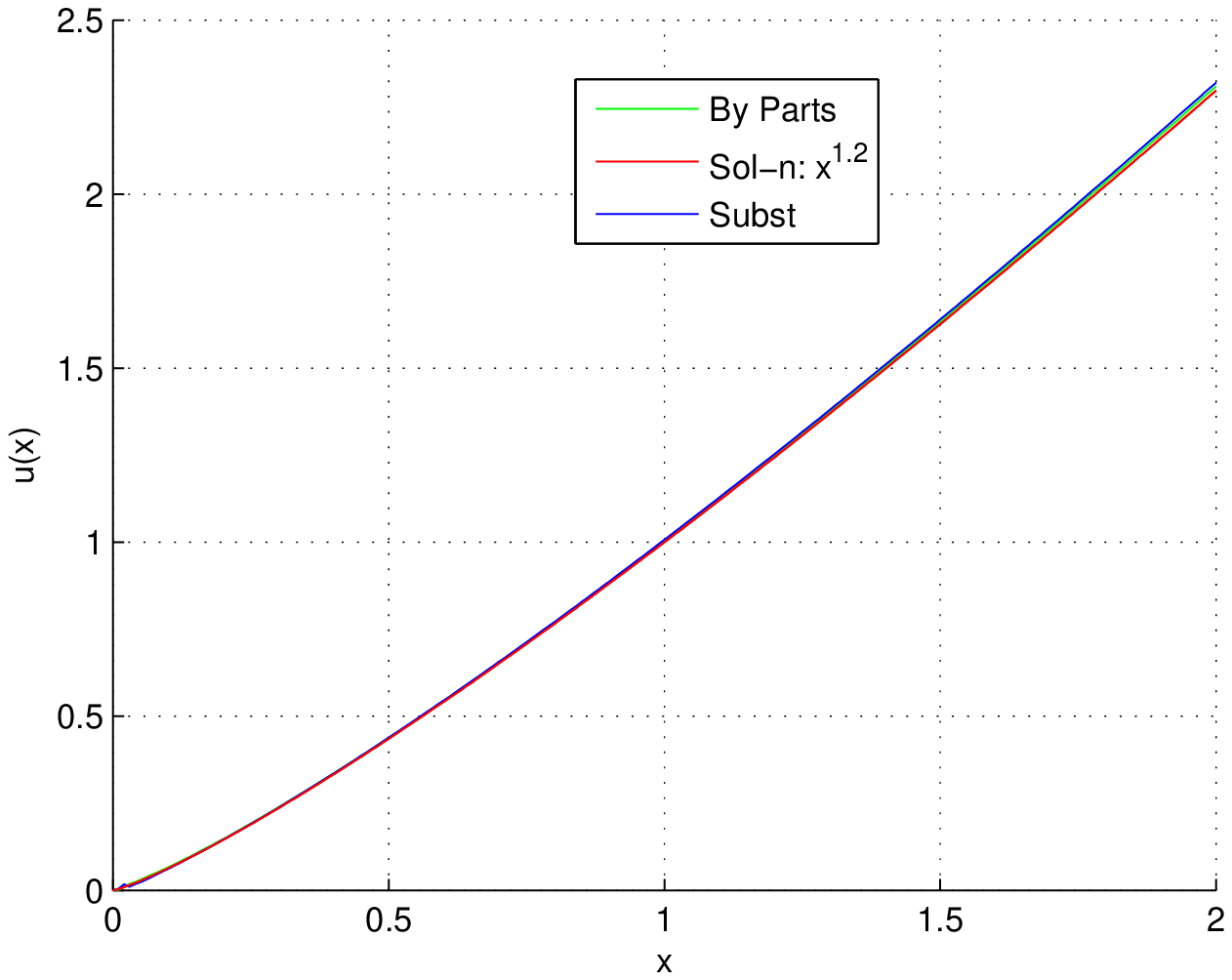}
		\includegraphics[width=7.5cm]{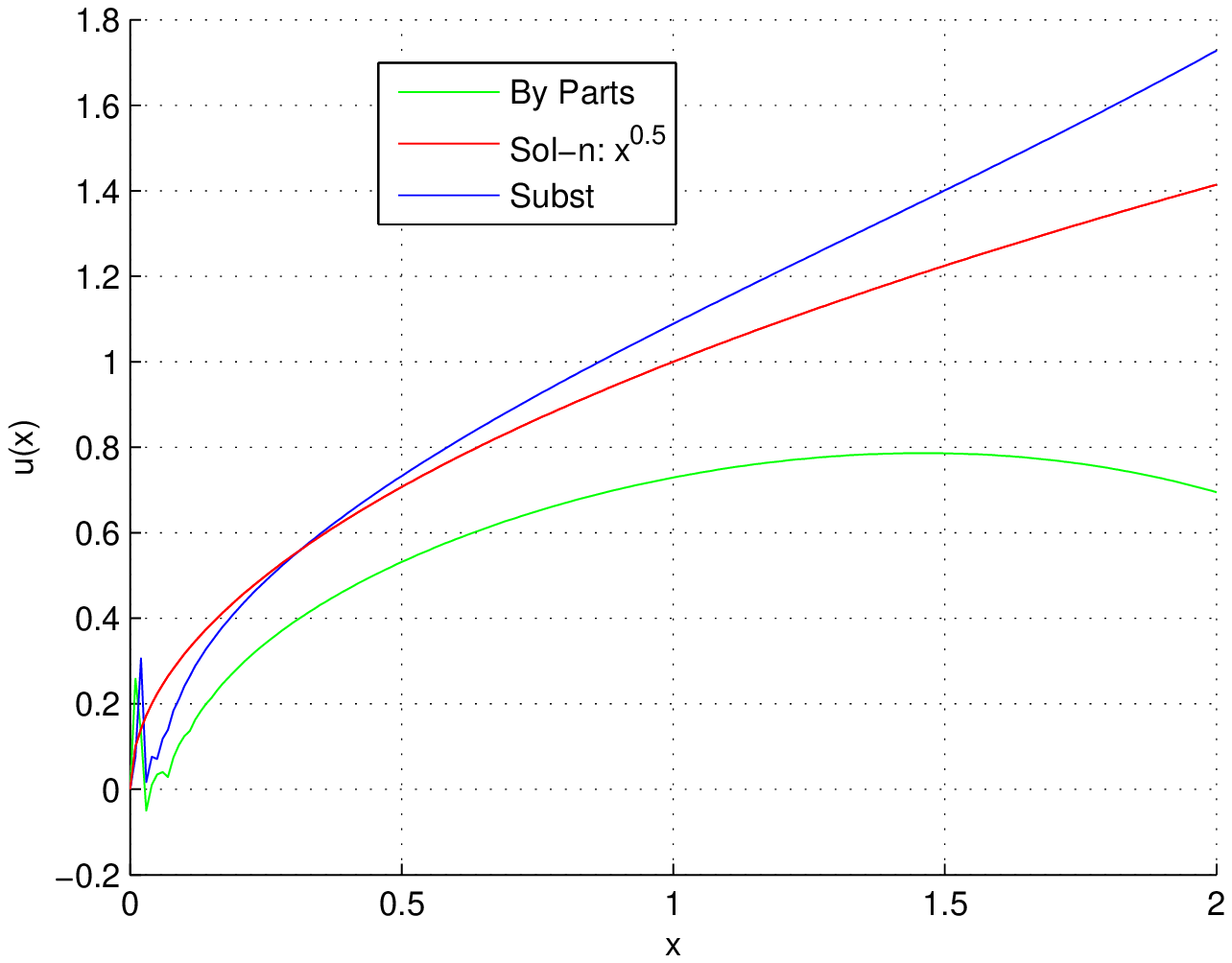}
		\caption{1a (left), equation (\ref{5-1a}): both substitution and integration by parts methods match each other and the exact solution.\\
		1b (right), equation (\ref{5-1b}): both substitution and by-parts methods diverge and don't match neither each other nor the exact solution.  \label{GoodEx}}	
	\end{center}
\end{figure}
\noindent In Figure 2 we consider slightly different equations 
\bq
(\frac{1}{4}+x^2)D^{0.5}u(x)-u(x)+x-\frac{2}{\sqrt \pi}(\frac{1}{4}\sqrt x +x^{2.5})=0 \label{5-2a}
\eq
with exact solution $u=x$. Also, we consider equation
\bq
(\frac{1}{100}+x^2)D^{0.5}u(x)-u(x)+x^2-\frac{8}{3\sqrt \pi}(\frac{1}{100}x^{1.5} +x^{3.5})=0. \label{5-2b}
\eq
with exact solution $u=x^2$.
In Figure 2a (equation \eqref{5-2a}) the method of substitution fails whereas the by-parts method is valid.
In Figure 2b (equation \eqref{5-2b}) we observe the opposite case: the method of substitution is valid but the method of integration by parts fails.

In all examples the computational step is $h=0.01$. As we can see, even though the equations look quite similar, the results of calculations differ: from almost perfect match by both methods to complete mismatch, which demonstrates that only if both methods produce almost identical results, then the computation is reliable. This is the additional justification of the proposed dual approach.

\begin{figure}[H] 
	\begin{center}
		\includegraphics[width=7.5cm]{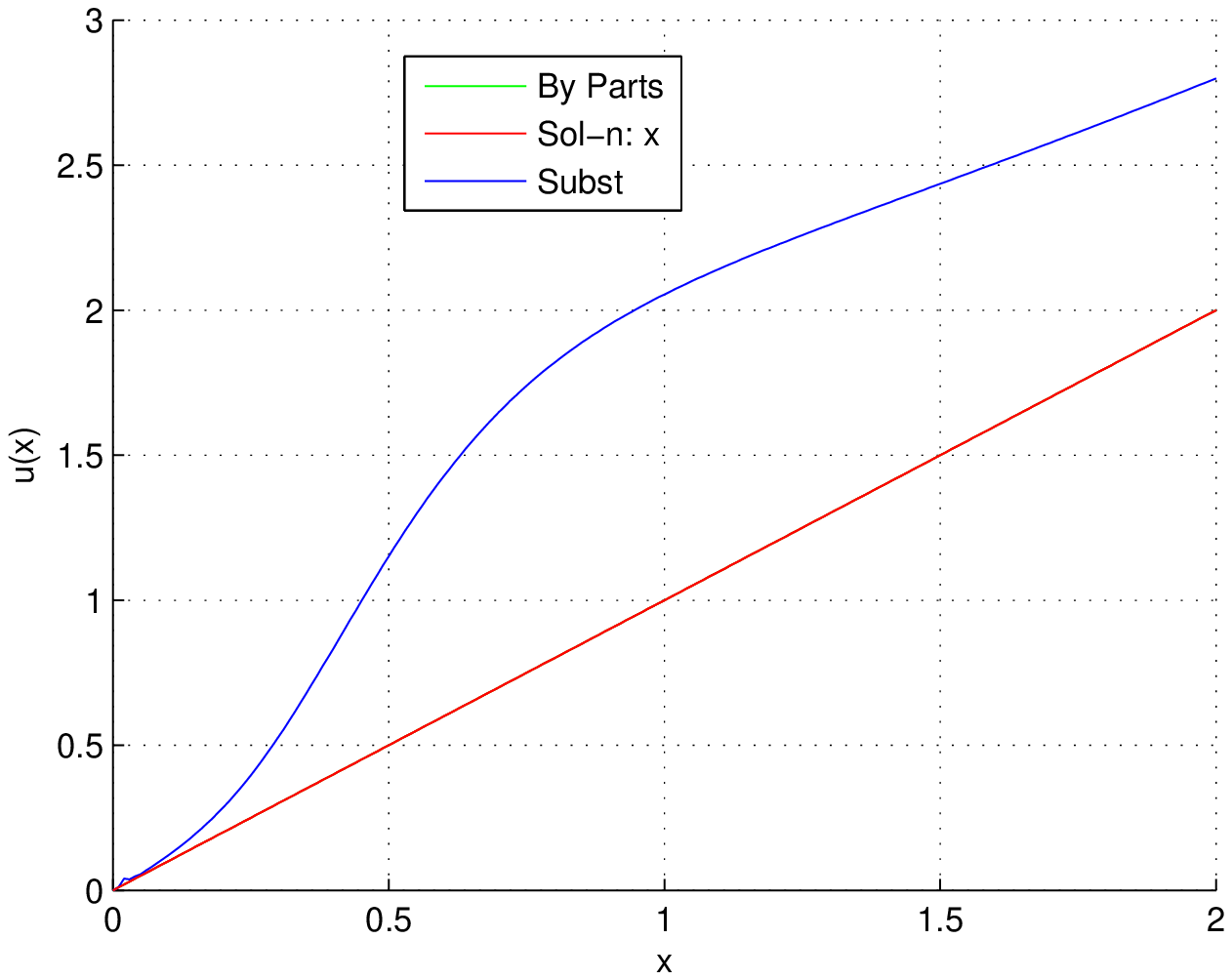}
		\includegraphics[width=7.5cm]{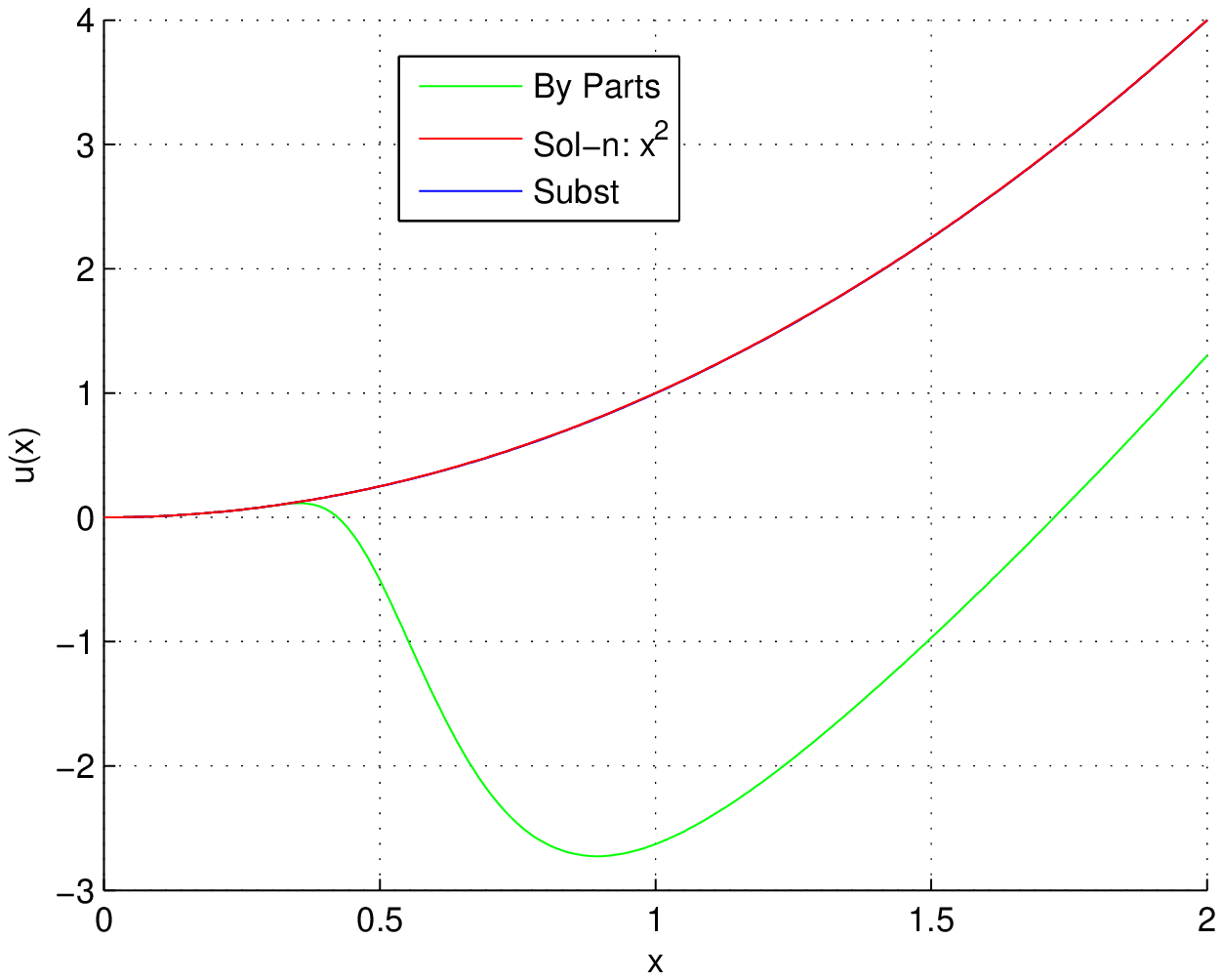}
		\caption{2a (left), equation (\ref{5-2a}): by-parts method matches the solution, but substitution fails.\\
		2b (right), equation (\ref{5-2b}): by-parts method fails, but the substitution method is valid.
		\label{ByPartsEx}}
	\end{center} 
\end{figure}

\section{Reliability for quasilinear equations}\label{SectExamples} 

\begin{figure}[H] 
	\begin{center}
		\includegraphics[width=12cm]{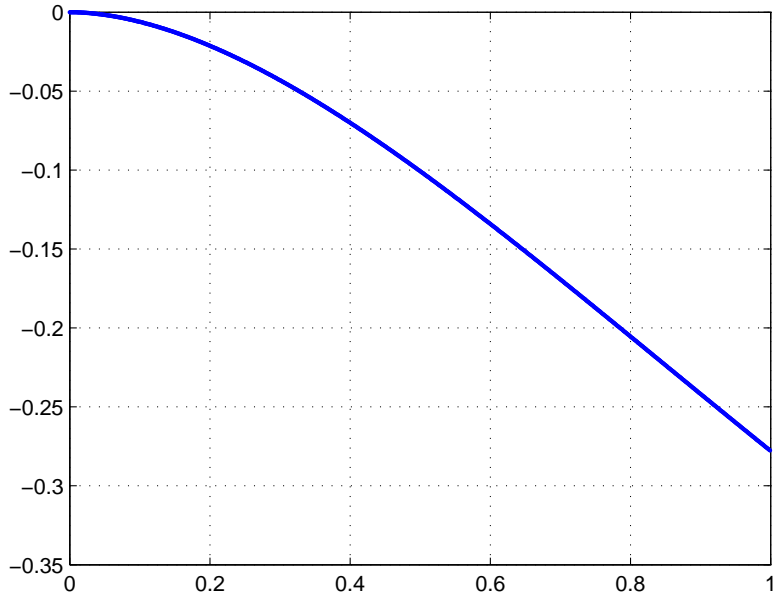}
\caption{Solution to equation \eqref{Example1} 
with step $h=0.001$ obtained by both methods, which provide almost identical results. The solution is reliable.}\label{AllAlphas}
	\end{center}
\end{figure}

{\bf Example 1.} Let's consider equation
 \begin{equation}\label{Example1}
 D^{0.3}u(x)+xD^{0.7}u(x)+\cos(u(x))D^{0.9}u(x)+\sin(x)=u^2(x)+\tan(u(x))
 \end{equation}
 with initial condition $u(0)=0$.  The calculated solution is represented in Figure \ref{AllAlphas}.
\noindent Table \ref{SolEq1} shows the results of solving equation \eqref{Example1} by two methods. The residual is calculated as the difference between the right and left sides of the equation. 

\begin{table}[H]
	\begin{center}
		\begin{tabular}{ | c | c | c || c | c |}\hline
    & Solution found & Residual           &  Solution found     & Residual   	\\
$x$ & using By-parts &                    &  using Substitution &                \\\hline
0.1 & -0.0061330982  & $3\times 10^{-9}$  &-0.0061330846        &  $2\times 10^{-7}$ \\\hline
0.2 & -0.0212821228  & $8\times 10^{-10}$ & -0.0212821387       & $1\times 10^{-7}$ \\\hline
0.3 & -0.0431124645  & $5\times 10^{-10}$ &-0.0431125027        & $1\times 10^{-7}$ \\\hline
0.4 & -0.0700231242  & $3\times 10^{-10}$ & -0.0700231812       & $1\times 10^{-7}$\\\hline
0.5 & -0.1007208712  & $3\times 10^{-10}$ & -0.1007209446       & $8\times 10^{-8}$ \\\hline
0.6 & -0.1341161406  & $2\times 10^{-10}$ & -0.1341162285       & $7\times 10^{-8}$ \\\hline
0.7 & -0.1692773586  & $2\times 10^{-10}$ & -0.1692774592       & $7\times 10^{-9}$ \\\hline
0.8 & -0.2054041267  & $2\times 10^{-10}$ & -0.2054042388       & $6\times 10^{-8}$ \\\hline
0.9 & -0.2418082833  & $2\times 10^{-10}$ & -0.2418084054       & $5\times 10^{-8}$\\\hline
1.0 & -0.2778991084  & $1\times 10^{-10}$ & -0.2778992392       & $8\times 10^{-8}$ \\\hline
		\end{tabular}
\caption{Solution to equation \eqref{Example1}.  Step $h=0.001$.  }\label{SolEq1}
\end{center}	
\end{table}
\noi As we can see in Table \ref{SolEq1}, the results of both methods are almost identical and residual $\sim 10^{-9}-10^{-7}$. However,  the true error is definitely bigger than $h^2=10^{-6}$, because on top of the error of the solution we must apply the approximation error for fractional derivatives. Consequently, there is no need for us to be delighted is the calculated residual is small. 
The actual error between the computations and the solution will be bigger than $O(h^2)$.\\

To support further the above warning on the computational errors, let us solve slightly changed equation \eqref{Example1} with the modified function $f(x)$ such that it produces the known result $u(x)=-x^2$. This value for 
$u(x)=-x^2$ is chosen because its graph is similar to the above computational results for \eqref{Example1}.  
Since $D^{0.3}(-x^2)=-\ds\frac{200}{119}x^{1.7}$, and $D^{0.7}(-x^2)=-\ds\frac{200}{39}x^{1.3}$ 
and $D^{0.9}(-x^2)=-\ds\frac{200}{11}x^{1.1}$, then we consider equation like \eqref{Example1}
\bqr
&& D^{0.3}u(x)+xD^{0.7}u(x)+\cos(u(x))D^{0.9}u(x)+x^4-\tan(x^2) \nonumber \\
&& +\frac{200}{119\Gamma(0.7)}x^{1.7}+\frac{200}{39\Gamma(0.3)}x^{2.3}+\cos(x^2)\frac{200}{11\Gamma(0.1)}x^{1.1}=
u^2(x)+\tan(u(x)) \label{Example1a}
\eqr
with the initial condition $u(0)=0$. Its exact solution is $u=-x^2$. In this case we can compare the calculated results by both methods with the known analytic solution:

\begin{table}[H]
	\begin{center}
		\begin{tabular}{ | c | c | c | c | c | c | c |}\hline
	  & Exact  & Solution found &               & Solution found  &       & Difference b/w \\
$x$ & $-x^2$ & using by-parts & Error         & by Substitution & Error & two methods \\\hline
0.1 & -0.01  & -0.0100507844  & $5.1 \times 10^{-5}$ & -0.0100508224   & $5.1 \times 10^{-5}$ & $3.8 \times 10^{-8}$ \\\hline
0.2 & -0.04 & -0.0400897392 & $9.0 \times 10^{-5}$ & -0.0400897712 & $9.0 \times 10^{-5}$ & $3.2 \times 10^{-8}$ \\\hline
0.3 & -0.09 & -0.0901232321 & $1.2 \times 10^{-4}$ & -0.0901232606 & $1.2 \times 10^{-4}$ & $2.9 \times 10^{-8}$ \\\hline
0.4 & -0.16 & -0.1601525607 & $1.5 \times 10^{-4}$ & -0.1601525867 & $1.5 \times 10^{-4}$ & $2.6 \times 10^{-8}$ \\\hline
0.5 & -0.25 & -0.2501785451 & $1.7 \times 10^{-4}$ & -0.2501785691 & $1.8 \times 10^{-4}$ & $2.4 \times 10^{-8}$ \\\hline
0.6 & -0.36 & -0.3602021066 & $2.0 \times 10^{-4}$ & -0.3602021290 & $2.0 \times 10^{-4}$ & $2.2 \times 10^{-8}$ \\\hline
0.7 & -0.49 & -0.4902245812 & $2.2 \times 10^{-4}$ & -0.4902246024 & $2.2 \times 10^{-4}$ & $2.1 \times 10^{-8}$ \\\hline
0.8 & -0.64 & -0.6402483126 & $2.5 \times 10^{-4}$ & -0.6402483332 & $2.5 \times 10^{-4}$ & $2.1 \times 10^{-8}$ \\\hline
0.9 & -0.81 & -0.8102786790 & $2.9 \times 10^{-4}$ & -0.8102786996 & $2.8 \times 10^{-4}$ & $2.1 \times 10^{-8}$ \\\hline
1.0 & -1.00 & -1.0003337914 & $3.3 \times 10^{-4}$ & -1.0003338137 & $3.3 \times 10^{-4}$ & $2.2 \times 10^{-8}$ \\\hline
		\end{tabular}
		\caption{Solution to equation \eqref{Example1a} found by each method with step $h=0.001$.}\label{SolEq1a}
	\end{center}	
\end{table}
As we can see in Table \ref{SolEq1a}, both methods produce almost identical results and, therefore, almost the same errors. The difference between the solutions is $\sim 10^{-8}$.  
However, the error is about $10^{-4}$, which is, as expected, more than the precision of calculation of fractional derivatives, but what's important is that the solutions 
found by both methods almost match. This indicates that it was found correctly and, consequently, is reliable.

Please note, that differential equations with $1<\alpha<2$ can be evaluated in the same fashion with similar precision.  For example, for equation
\begin{equation}\label{Example2}
xu(x) D^{1.3}u(x)+e^x D^{1.7}u(x)-40x^4=\sin(u(x))
\end{equation} 
 with initial conditions $u(0)=0, u'(0)=0$ both methods produce almost identical results.\\

\noindent {\bf Example 2.}
For equations containing fractional derivatives of different orders (between $0<\alpha<1$ and $1<\alpha<2$), 
the results may be less reliable.  Having two different methods of calculation provides us with a path to analyzing the validity of the results.

\begin{figure}[H] 
	\begin{center}
		\includegraphics[width=12cm]{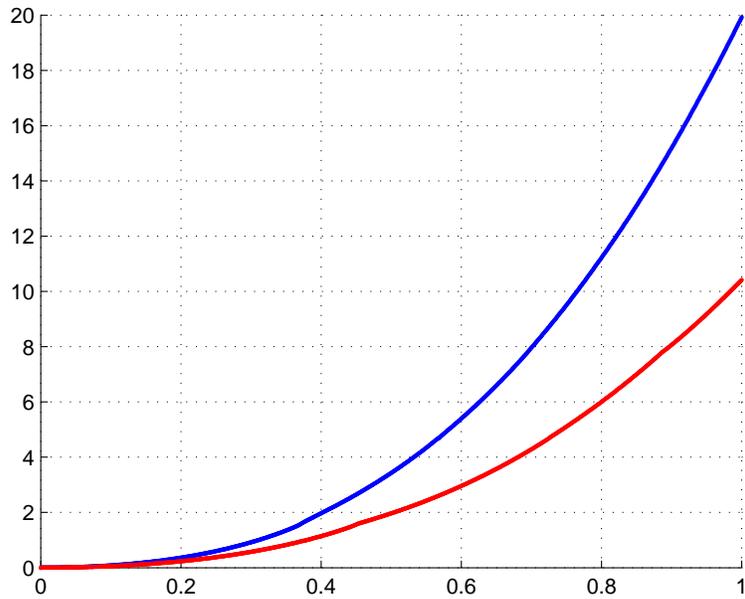}
		\caption{'Solutions' to equation \eqref{Example3bad} with step $h=0.001$.  Red line -- substitution method, blue -- by-parts method.  The solutions are not reliable.}\label{AllBad}
	\end{center}
\end{figure}

\begin{figure}[H] 
	\begin{center}
		\includegraphics[width=12cm]{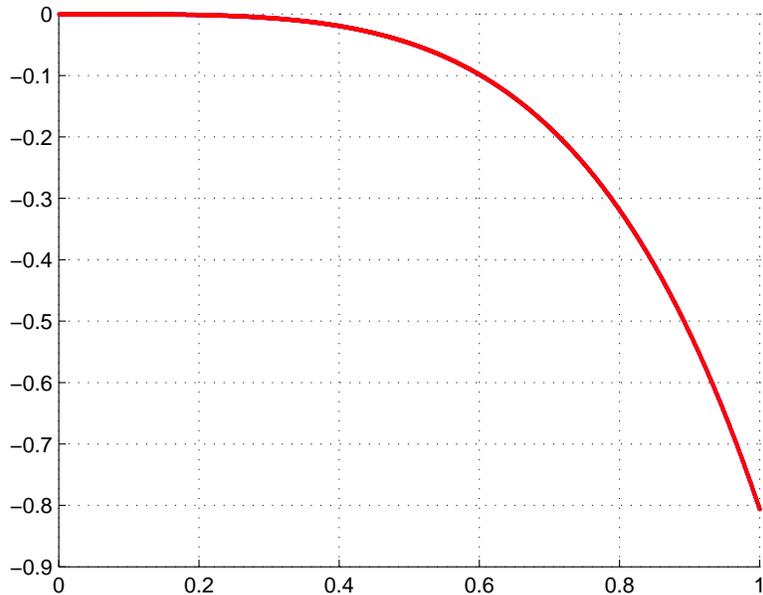}
		\caption{Solution for \eqref{Example3} 
		with step $h=0.001$. 
		The methods produce almost identical result, which is reliable.}\label{AllAlphas3}
	\end{center}
\end{figure}

If we solve semilinear equation
\begin{equation}\label{Example3bad}
x D^{1.7}u(x)+x^2 D^{0.3}u(x)-3x=5u(x)+\tan(u(x)),
\end{equation}
then our two methods produce quite different outputs, which indicate that the found by either method solution is not reliable (see Figure \ref{AllBad}).\\
If we slightly modify the entries in \eqref{Example3bad} the computation becomes more stable. For example, by changing $-3x$ to $3x^3+\tan(x^4)$, we arrive at equation
\begin{equation}\label{Example3}
	x D^{1.7}u(x) + x^2 D^{0.3}u(x)+3x^3+\tan(x^4)=5u(x)+\tan(u(x)).
\end{equation}
The solution for equation \eqref{Example3} is presented in Figure \ref{AllAlphas3} using both methods.
For further analysis, we replace terms $3x^3+\tan(x^4)$ in \eqref{Example3} by similar function, which provides exact solution $u(x)=x^3$. Thus, we consider the following equation:
\begin{equation}\label{Example3chk}
  x D^{1.7}u(x) + x^2 D^{0.3}u(x)+ 5x^3+\tan(x^3)-\frac{2000x^{4.7}}{1071\Gamma(0.7)}-\frac{200x^{2.3}}{13\Gamma(0.3)}=5u(x)+\tan(u(x)).
\end{equation}
Numerical solution of equation \eqref{Example3chk} produces $u(x)=x^3$ with error $\sim 10^{-4}$ by both almost coinciding methods.  This confirms that solution to \eqref{Example3} is reliable, even though close equation \eqref{Example3bad} suggests opposite. These examples justify once again the necessity of the proposed dual approach. 

\section{Conclusion}
To analyze the computational methods for fractional differential equations, we introduce the substitution numerical method and, along with well-known by-parts method, suggest the dual approach for the reliability of computations for linear and quasilinear fractional differential equations. We show numerically that in the cases when both methods produce almost identical outputs, the solution is reliable, and its level of accuracy is close to the accuracy of the approximations of fractional derivatives. As a tool to verify the accuracy for fractional derivatives, we use the Taylor series expansion. 
%
We demonstrate the validity of the dual approach by solving linear and nonlinear fractional differential equations.
Our computational examples and counterexamples demonstrate the necessity of the dual approach.

\nocite{*}


\begin{thebibliography}{9}
	\bibitem{Ramzi}
		R. Albadarneh, M. Zerqat, I. Batiha.
		\textit {Numerical Solutions for linear and non-linear Fractional Differential Equations}.
		International Journal of Pure and Applied Mathematics, 106 No. 3, 2016, pp. 859--871.

	\bibitem{Ramzi2}
		R. Albadarneh, M. Zurigat, I. Batiha.
		\textit{Numerical Solutions for linear Fractional Differential Equations of order $1<\alpha < 2$ using Finite Difference Method (FFDM)}.
		J. Math. Computer Sci. 16 (2016), 103--111.

\bibitem{Boyadzhiev}
S.P. Mirevski, L. Boyadjiev, R. Scherer.
	\textit{On the Riemann-Louisville	fractional calculus, g-Jacobi functions and F-Gauss functions}. 
	Applied	Mathematics and Computation 187, No 1, 2007, pp.315–325
	
\bibitem{Gerasimov} A.N.Gerasimov. \textit{A generalization of the deformation laws and its application to the problems of internal friction}. Applied Mathematics and Mechanics 12, 1948, pp. 251--260.
	
	
\bibitem{Hilfer} 	Z. Tomovski, R. Hilfer, H.M. Srivastava. 
	\textit{Fractional and operational calculus with generalized fractional derivative operators and Mittag–Leffler type functions}.
	Integral Transforms and Special Functions, 2010, 21(11), pp.797-814

	\bibitem{Kilbas}
	A.A. Kilbas, H.M. Srivastava, J.J. Trujillo.
	\textit{Theory and applications of Fractional Differential equations}.
	Fakulteit der Exacte Wetenschappen, 2006, Amsterdam, The Netherlands.
	

	\bibitem {Mathews}
		J. Mathews, K. Fink.
		\textit{Numerical Methods Using MATLAB}.
		Fourth Edition.

	\bibitem{Podlubny}
		I. Podlubny.
		\textit{Fractional Differential Equations}.
		Mathematics in Science and Engineering, Vol. 198, 1999.

\bibitem{Zhuk-Kilbas} N.V.Zhukovskaya, A.A.Kilbas. \textit{Soilving homogeneous fractional differential equations of Euler type}. Differential Equations 47, No. 12, 2011, pp. 1714--1725. 


\end{thebibliography}
\end{document}